\documentclass[12pt]{amsart}
\usepackage{latexsym}
\usepackage{amssymb}
\usepackage{amsmath}
\usepackage{mathrsfs}
\usepackage{bbm}
\usepackage{bbold}
\def\epi{\mathop{\fam0 epi}\nolimits}
\def\cop{\mathop{\fam0 cop}\nolimits}
\def\dom{\mathop{\fam0 dom}\nolimits}

\begin{document}

\title{
Abstract Convexity and \hbox{Cone-Vexing Abstractions}
}

\author{S.~S. Kutateladze}
\address[]{
Sobolev Institute of Mathematics\newline
\indent 4 Koptyug Avenue\newline
\indent Novosibirsk, 630090\newline
\indent Russia}
\email{
sskut@member.ams.org
}
\date{May 18, 2007}
\thanks{I am grateful to Professor G.~Litvinov who kindly
invited this talk to the ``tropics.''}
\maketitle

This talk is devoted to some origins of abstract convexity and
a~few vexing limitations on the range of abstraction in convexity.
Convexity is a relatively recent subject. Although the noble
objects of Euclidean geometry are mostly convex, the abstract
notion of a convex set appears only after the Cantor paradise was
founded. The idea of convexity feeds generation, separation, calculus, and
approximation.  Generation appears as duality; separation, as optimality;
calculus,  as representation; and approximation, as stability.

{\bf 1. Generation.}  Let $\overline{E}$  be a~complete lattice
 $E$ with the adjoint top $\top:=+\infty$ and bottom $\bot:=-\infty$.
Unless otherwise stated, $Y$ is usually a~{\it Kantorovich space\/}
which is a Dedekind complete vector lattice in another terminology.
Assume further that $H$  is some subset of $E$ which is by implication a~(convex)
cone in $E$, and so the bottom of $E$
lies beyond~$H$. A subset $U$  of~ $H$ is {\it convex relative to~}
$H$ or $H$-{\it convex\/}, in symbols $U\in\mathfrak{V}(H,\overline{E})$,
provided that $U$ is the $H$-{\it support set\/}
$U^H_p:=\{h\in H:h\le p\}$ of some element $p$ of $\overline{E}$.

Alongside the $H$-convex sets we consider
the so-called $H$-convex elements. An element   $p\in \overline{E}$
is  $H$-{\it convex} provided that $p=\sup U^H_p$;~i.e., $p$
represents the supremum of the $H$-support set of~$p$.
The $H$-convex elements comprise the cone which is denoted by
$\mathscr C(H,\overline{E}$).  We may omit  the references to $H$ when $H$ is clear
from the context. It is worth noting that
convex elements and sets are ``glued together''
by the {\it Minkowski diality\/} $ \varphi:p\mapsto U^H_p$.
This duality enables us to study convex elements and sets simultaneously.

Since the classical results by Fenchel \cite{Fenchel}
and H\"ormander
\cite{Her, Notions} it has been well known
that the most convenient and conventional classes of convex  functions
and sets are  $\mathscr C(A(X),\overline{\mathbb R^X} )$ and $\mathfrak{V}(X',\overline{\mathbb R^X})$.
Here    $X$ is a locally convex space,   $X'$ is the dual of~$X$,
and $A(X)$ is the space of affine functions on  $X$
(isomorphic with  $X'\times \mathbb R$).

In the first case the Minkowski duality is the mapping
 $f\mapsto\text{epi} (f^*)$ where
$$
f^*(y):=\sup\limits_{x\in X}(\langle y,x\rangle - f(x))
$$
is the {\it Young--Fenchel transform\/} of~$f$ or the {\it conjugate function\/} of~$f$~.
In the second case we prefer to write down the inverse of the Minkowski
duality which sends  $U$  in $\mathfrak{V}(X',\overline{\mathbb R}^X)$
to the standard {\it support function}
$$
\varphi^{-1}(U):x\mapsto\sup\limits_{y\in U}
\langle y,x\rangle.
$$
As usual, $\langle\cdot,\cdot\rangle$ stands for the canonical pairing
of~$X'$ and~$X$.

This idea of abstract convexity lies behind many current objects
of analysis and geometry. Among them we list the ``economical'' sets
with boundary points meeting the Pareto criterion, capacities, monotone
seminorms, various classes of functions convex in some generalized sense,
for instance, the Bauer convexity in Choquet  theory, etc.
It is curious that there are ordered vector spaces consisting of
the convex elements with respect to  narrow cones with finite generators.
Abstract convexity is traced and reflected, for instance, in
\cite{MD}--%
\cite{IoRu}.

{\bf 2. Separation.}
Consider cones
$K_1$
and
$K_2$
in a topological vector space
$X$
and put
$\varkappa:=(K_1,K_2)$.
Given a pair
$\varkappa$
define the correspondence
$\Phi_{\varkappa}$
from
$X^2$
into
$X$
by the formula
$$
\Phi_{\varkappa}:=\{(k_1,k_2,x)\in X^3:
x=k_1-k_2\in K_\imath\ (\imath:=1,2)\}.
$$
Clearly, $\Phi_{\varkappa}$ is a cone or, in other words,
a~conic correspondence.

The pair $\varkappa$ is {\it nonoblate\/}
whenever $\Phi_{\varkappa}$
is open at the zero. Since
$\Phi_{\varkappa}(V)=V\cap K_1-V\cap K_2$
for every
$V\subset X$,
the nonoblateness of $\varkappa$
means that $$
\varkappa V:=(V\cap K_1-V\cap K_2)\cap(V\cap K_2-V\cap K_1)
$$
is a zero neighborhood  for every zero neighborhood~
$V\subset X$.
Since $\varkappa V\subset V-V$, the nonoblateness of
$\varkappa$ is equivalent to the fact that the system of sets
$\{\varkappa V\}$ serves as a filterbase of zero neighborhoods while
$V$ ranges over some base of the same filter.

Let $\Delta_n:x\mapsto(x,\dots,x)$ be
 the embedding of
$X$
into the diagonal
$\Delta_n(X)$
of
$X^n$.
A pair of cones
$\varkappa:=(K_1,K_2)$
is nonoblate if and only if
$\lambda:=(K_1\times K_2,\Delta_2(X))$
is nonoblate in~$X^2$.

Cones
$K_1$
and
$K_2$
constitute a nonoblate pair if and only if the conic correspondence~
$\Phi\subset X\times X^2$
defined as
$$
\Phi:=\{(h,x_1,x_2)\in X\times X^2 :
x_\imath+h\in K_\imath\ (\imath:=1,2)\}
$$
is open at the zero. Recall that a convex correspondence
$\Phi$ from $X$
into $Y$ is open at the zero if and only if the H\"ormander transform
of $X\times\Phi$ and the cone
$\Delta_2(X)\times\{0\}\times\mathbb R^+$
constitute a nonoblate pair in~$X^2\times Y\times\mathbb R$.

Cones~
$K_1$
and~
$K_2$
in a topological vector space~
$X$
are {\it in general position\/}
provided that

{\bf (1)}~
the algebraic span of $K_1$
and~
$K_2$
is some subspace
$X_0\subset X$;
i.e.,
$X_0=K_1-K_2=K_2-K_1$;

{\bf (2)}~the subspace
$X_0$
is complemented; i.e., there exists a continuous projection
$P:X\rightarrow X$
such that
$P(X)=X_0$;

{\bf (3)}~$K_1$
and~
$K_2$
constitute a nonoblate pair in~
$X_0$.

Let $\sigma_n$
stand for the rearrangement of coordinates
$$
\sigma_n:((x_1,y_1),\dots, (x_n,y_n))\mapsto ((x_1,\dots,x_n),
(y_1,\dots,y_n))
$$
which establishes an isomorphism between
$(X\times Y)^n$
and
$X^n\times Y^n$.

Sublinear operators  $P_1,\dots,P_n:X\rightarrow E\cup \{+\infty\}$
are {\it in general position\/} if so are
the cones $\Delta_n(X)\times E^n$
and
$\sigma_n(\epi (P_1)\times\dots\times\epi (P_n))$.
A similar terminology applies to convex operators.

Given a cone $K\subset X$, put
$$
\pi_E(K):=\{T\in\mathscr L(X,E): Tk\leq 0\ (k\in K)\}.
$$
We readily see that
$\pi_E(K)$
is a cone in
$\mathscr L(X,E)$.

{\scshape Theorem.} {\sl Let
$K_1,\dots,K_n$
be cones in a topological vector space~
$X$
and let
$E$
be a topological Kantorovich space.  If
$K_1,\dots,K_n$
are in general position then}
$$
\pi_E(K_1\cap\dots\cap K_n)=\pi_E(K_1)+\dots+\pi_E(K_n).
$$
This formula opens  a way to various separation results.

{\scshape Sandwich Theorem.} {\sl Let
$P,Q:X\rightarrow E\cup \{+\infty\} $
be sublinear operators in general position.
If
$P(x)+Q(x)\geq 0$
for all
$x\in X$
then there exists a~continuous linear operator~
$T:X\rightarrow E$
such that}
$$
-Q(x)\leq Tx\leq P(x)\quad  (x\in X).
$$

Many efforts were made to abstract these results to
a more general algebraic setting and, primarily,
to semigroups. The relevant separation results are
collected  in~\cite{Fuch}.

{\bf 3. Calculus.}
Consider a~Kantorovich space
$E$
and an arbitrary nonempty set
$\mathfrak A$.
Denote by
$l_\infty (\mathfrak A,E)$
the set of all order bounded mappings from
$\mathfrak A$
into $E$; i.e.,
$f\in l_\infty (\mathfrak A,E)$
if and only if
$f:\mathfrak A \to E$
and the set
$\{f(\alpha):\alpha\in\mathfrak A\}$
is order bounded in $E$.
It is easy to verify that
$l_\infty (\mathfrak A,E)$ becomes a Kantorovich
space if endowed with the coordinatewise algebraic operations and order.
The operator
$\varepsilon_{\mathfrak A, E}$
acting from
$l_\infty (\mathfrak A,E)$
into
$E$
by the rule
$$
\varepsilon_{\mathfrak A, E}:f\mapsto\sup \{f(\alpha):\alpha\in\mathfrak A\}
\quad (f\in l_\infty (\mathfrak A,E))
$$
is called the {\it canonical sublinear operator\/}
given
$\mathfrak A$
and
$E$.
We often write
$\varepsilon_{\mathfrak A}$
instead of
$\varepsilon_{\mathfrak A, E}$
when it is clear from the context what Kantorovich space is meant.
The notation
$\varepsilon_n$
is used when the cardinality of~$\mathfrak A$
equals
$n$ and we call the operator
$\varepsilon_n$
{\it finitely-generated}.

Let
$X$
and
$E$
be ordered vector spaces. An operator
$p: X\to E$
is called {\it increasing\/} or
{\it isotonic\/}
 if for all
$x_1, x_2 \in X$
from
$x_1\leq x_2$
it follows
that
$p(x_1)\leq p(x_2)$.
An increasing linear operator is also called {\it positive}.
As usual, the collection of all positive linear operators in the space
$L(X,E)$ of all linear operators is denoted
by
$L^+ (X,E)$.
Obviously, the positivity of a~linear operator
$T$ amounts to the
inclusion
$T(X^+) \subset E^+$,
where
$X^+:=\{x\in X: x\geq 0\}$
and
$E^+:=\{e\in E: e\geq 0\}$
are the
{\it positive cones\/}
in
$X$
and
$E$
respectively.
Observe that every canonical operator is increasing and sublinear,
while every finitely-generated canonical operator is order continuous.

Recall that
$\partial p:=\partial p(0)=\{ T \in L (X,E):$ $(\forall x
\in X  )\  T x\leq p(x)\}$
is the {\it subdifferential\/}
at the zero
or  {\it support
set\/}
of a~sublinear operator~
$p$.

Consider a~set~
$\mathfrak A$
of linear operators acting from
a~vector space
$X$
into a~Kantorovich space
$E$.
The set
$\mathfrak A$
is {\it weakly order  bounded\/} if
the set
$\{\alpha x:\alpha\in\mathfrak A\}$
is order bounded for every
$x\in X$.
We denote by
$\langle\mathfrak A\rangle x$
the mapping that assigns the element
$\alpha x\in E$
to each
$\alpha\in \mathfrak A$,
i.e.
$\langle\mathfrak A\rangle x: \alpha\mapsto\alpha x$.
If
$\mathfrak A$
is weakly order bounded then
$\langle\mathfrak A\rangle x\in l_\infty (\mathfrak A,E)$
for every fixed
$x\in X$.
Consequently, we obtain the linear operator
$\langle\mathfrak A\rangle:X\rightarrow l_\infty (\mathfrak A,E)$
that acts as
$\langle\mathfrak A\rangle:x\mapsto\langle\mathfrak A\rangle x$.
Associate with
$\mathfrak A$
one more operator
$$
p_{\mathfrak A}: x\mapsto\sup \{\alpha x: \alpha\in\mathfrak A\}\quad(x\in X).
$$
The operator
$p_{\mathfrak A}$
is sublinear. The support set
$\partial p_{\mathfrak A}$
is denoted by
$\cop (\mathfrak A)$
and referred to as  the {\it support hull\/} of
$\mathfrak A$.
These definitions  entail the following

{\scshape Theorem.} {\sl If
$p$
is a~sublinear operator with
$\partial p=\cop (\mathfrak A)$
then $
P=\varepsilon_{\mathfrak A}\circ \langle\mathfrak A\rangle.
$
Assume further that
$p_1: X\to E$
is a~sublinear operator and
$p_2: E\to F$
is an increasing sublinear operator. Then
$$
\partial (p_2\circ p_1)=\left\{ T\circ\langle\partial p_1\rangle: T\in L^+
(l_{\infty}(\partial p_1, E),F)\ \wedge\ T\circ\Delta_{\partial p_1}\in
\partial p_2 \right\}.
$$
Furthermore, if
$\partial p_1=\cop (\mathfrak A_1)$
and
$\partial p_2=\cop (\mathfrak A_2)$
then}
$$
\gathered
\partial (p_2\circ p_1)
=\bigl\{T\circ\langle\mathfrak A_1\rangle : T\in L^+
(l_{\infty}(\mathfrak A_1,E),F)\
\\
\wedge\
\left(\exists\alpha\in\partial\varepsilon_{\mathfrak A_2}\bigr)\
T\circ\Delta_{\mathfrak A_1}=\alpha\circ\langle\mathfrak A_2\rangle\right\}.
\endgathered
$$

More details on subdifferential calculus and applications to optimality
are collected in~\cite{Subdiff}.

{\bf 4. Approximation.} Study of stability in abstract convexity
is accomplished sometimes by introducing various  epsilons
in appropriate places. One of the earliest attempts in this direction
is connected with the classical  Hyers--Ulam stability theorem for
$\varepsilon$-convex functions. The most recent results are collected
in~\cite{Almost}. Exact calculations with epsilons and sharp estimates
are sometimes bulky and slightly mysterious.  Some alternatives are suggested
by actual infinities, which is illustrated with the conception
of {\it infinitesimal optimality}.

Assume given a ~convex operator
$f:X\to E\cup{+\infty}$
and a~ point
$\overline x$
in the effective domain
$\dom(f):=\{x\in X:f(x)<+\infty\}$
of
~$f$.
Given
$\varepsilon \ge 0$
in the positive cone
$E_+$
of
$E$,
by the
$\varepsilon $-{\it subdifferential\/}
of~$f$
at
~$\overline x$
we mean the set
$$
\partial\, {}^\varepsilon\!f(\overline x):=\big\{T\in L(X,E):
(\forall x\in X)(Tx-Fx\le T\overline x -f\overline x+\varepsilon) \big\},
$$
with
$L(X,E)$
standing as usual for the space of linear operators
from~
$X$
to
~$E$.

Distinguish
some downward-filtered subset
~$\mathscr E$ of
$E$
that is composed of positive elements.
Assuming
$E$ and~$\mathscr E$
 standard, define the {\it monad\/}
$\mu (\mathscr E)$ of $\mathscr E$ as
$\mu (\mathscr E):=\bigcap\{[0,\varepsilon ]:\varepsilon \in
{}^\circ\!\mathscr E\}$.
The members of $\mu(\mathscr E)$ are {\it positive
infinitesimals\/}  with respect to~$\mathscr E $.
As usual,
${}^\circ\!\mathscr E$
denotes the external set of~
all standard members of
~$E$,
the {\it standard part\/} of
~$\mathscr E$.

We will agree that the monad $\mu (\mathscr E )$
is an external cone over ${}^\circ \mathbb R $ and, moreover,
$\mu (\mathscr E)\cap{}^\circ\! E=0$.
In application, $\mathscr E $ is usually
the filter of order-units of $E$.
The relation of
{\it infinite proximity\/} or
{\it infinite closeness\/}
between the members of $E$ is introduced as follows:
$$
e_1 \approx e_2 \leftrightarrow e_1 -e_2 \in\mu
(\mathscr E )\wedge e_2 -e_1 \in\mu (\mathscr E ).
$$

Since
$$
\bigcap\limits_{\varepsilon \in{}^\circ \mathscr E }\,
\partial _\varepsilon f(\overline x)=
\bigcup\limits_{\varepsilon \in\mu (\mathscr E )}\,
\partial _\varepsilon f(\overline x);
$$
therefore, the external set on both sides  is
 the so-called {\it infinitesimal subdifferential} of
$f$ at  $\overline x$. We denote this set by
$Df(\overline x)$.
The elements of
$Df(\overline x)$ are
{\it infinitesimal subgradients}
of $f$ at
~$\overline x$.
If the zero oiperator is an infinitesimal subgradient
of $f$ at $\overline x$ then $\overline x$ is called
an {\it infinitesimal minimum point\/} of $f$.
We abstain from indicating  $\mathscr E$ explicitly
since this leads to no confusion.

{\scshape Theorem.} {\sl
Let $f_1:X\times Y\rightarrow E\cup +\infty$ and
$f_2:Y\times Z\rightarrow E\cup +\infty$ be convex operators.
Suppose that the convolution
$f_2\vartriangle f_1$ is infinitesimally exact at some point $(x,y,z)$; i.e.,
$
(f_2\vartriangle f_1)(x,y)\approx f_1(x,y)+f_2(y,z).
$
If, moreover, the
convex sets $\epi(f_1,Z)$ and $\epi(X,f_2)$ are
in general
position then}
$$
D(f_2\vartriangle f_1)(x,y)=
Df_2(y,z)\circ Df_1(x,y).
$$


\bibliographystyle{plain}

\end{document}